\documentclass[11pt,oneside,reqno]{amsart}
\usepackage{bm,cancel,ulem}
\usepackage{amssymb,amsfonts}
\usepackage{amsmath,latexsym,enumerate}
\usepackage{mathrsfs}

\usepackage[top=1in, bottom=1in, left=1.25in, right=1.25in, marginparwidth=1.1in, marginparsep=0.05in]{geometry}

\usepackage[unicode,breaklinks=true,colorlinks=true]{hyperref}

\usepackage[svgnames]{xcolor}

\newtheorem{thm}{Theorem}[section]
\newtheorem{cor}[thm]{Corollary}
\newtheorem{lem}[thm]{Lemma}

\theoremstyle{definition}
\newtheorem{defn}[thm]{Definition}
\newtheorem{rem}[thm]{Remark}
\numberwithin{equation}{section}

\newcommand{\R}{\mathbb R}

\newcommand\NN{\mathbb{N}}

\newcommand{\cF}{\mathcal{F}}

\newcommand{\dy}{\,{\rm d}y}
\newcommand{\dx}{\,{\rm d}x}

\newcommand{\ds}{\,{\rm d}s}

\newcommand{\dtau}{\,{\rm d}\tau}
\newcommand{\dxi}{\,{\rm d}\xi}

\usepackage{accents}
\newcommand*{\dotc}[1] {\accentset{\mbox{\large\bfseries .}}{#1}}

\newcommand{\lec}{\lesssim}

\newcommand{\norm}[1]{\left\|#1\right\|}

\newcommand{\pd}{\partial}
\newcommand{\Om}{\Omega}

\newcommand{\tsum}{{\textstyle \sum}}
\renewcommand\rightmark{Maximal regularity for parabolic boundary value problems}

\begin{document}
\title{Remarks on the maximal regularity for parabolic boundary value problems with inhomogeneous data}%
\author[H. Chen]{Hui Chen}%
\address[H. Chen]
{School of Science, Zhejiang University of Science and Technology, Hangzhou, 310023, China }
\email{chenhui@zust.edu.cn}
\author[S. Liang]{Su Liang}
\address[S. Liang]
{Department of Mathematics, University of British Columbia, Vancouver, BC V6T1Z2, Canada }
\email{liangsu96@math.ubc.ca}
\author[T.-P. Tsai]{Tai-Peng Tsai}
\address[T.-P. Tsai]
{Department of Mathematics, University of British Columbia, Vancouver, BC V6T1Z2, Canada }
\email{ttsai@math.ubc.ca}

\date{}  

\begin{abstract}
Inspired by Ogawa-Shimizu \cite{Ogawa-Shimizu-2022} and Chen-Liang-Tsai \cite{Chen-Liang-Tsai-2025} on the second and first order derivative estimates of solution of heat equation in the upper half space with boundary data in homogeneous Besov spaces, we extend the estimates to any order of derivatives, including fractional derivatives.
\end{abstract}

\maketitle

\noindent {{\sl Key words:} Parabolic boundary value problems, maximal regularity, anisotropic Besov spaces}


\noindent {\sl AMS Subject Classification (2020):}  35K20, 46E35, 46B25

\tableofcontents

\section{Introduction}

For $d\geq 1$, consider the heat equation in the upper half space $\R^d_+$ with Dirichlet boundary condition,
\begin{equation}\label{hbvp}
\left\{
\begin{aligned}
\,\pd_t v-\Delta v&=0, \quad \hbox{in }\R^d_+\times \R,\\[4pt]
 v&=g, \quad \hbox{on }\pd \R^d_+\times \R.
\end{aligned}
\right.
\end{equation}
Here $\R^d_+=\{ x=(x',x_d): x'=(x_1,\ldots,x_{d-1})\in \R^{d-1},\ x_d>0\}$ when $d \ge 2$ and $\R^1_+ = (0,\infty)$.
Let $\Gamma$ be the heat kernel in $\R^d$, i.e.,
$\Gamma(x,t)=(4\pi t)^{-d/2} e^{-|x|^2/4t}$ for $t>0$ and $\Gamma(x,t)=0$ for $t \le 0$.
A solution of \eqref{hbvp} can be expressed in terms of its boundary data by
\begin{align}\label{def-v}
	v(x,t)=-2\int_\R\int_{\R^{d-1}}\pd_d \Gamma(x'-y',x_d,t-s)g(y',s)\dy'\ds.
\end{align}
When $d=1$, the variable $y'$ is null and \eqref{def-v} is reduced to $v(x,t)=-2\int_\R \pd_x \Gamma(x,t-s)g(s)\ds$.

In this note, we are concerned with the maximal regularity of parabolic boundary value problems. 
The maximal regularity of parabolic problems in the whole space and in domains has a long history. For parabolic maximal regularity with mixed exponents and nonzero boundary values, see e.g.~\cite{DHP-07,Ogawa-Shimizu-2022,Weidemaier-05}.

In particular, we establish the following estimates with boundary data $g$ in homogeneous anisotropic Besov spaces (see Section \ref{sec2}).
Let $\NN_0=\{0,1,2,\ldots\}$.

\begin{thm}\label{thm_main}
Suppose $v$ satisfies \eqref{def-v}, which solves heat equation with Dirichlet data $g$.
\begin{enumerate}[(i)]
\item For 
$1\leq p<\infty $, any $\alpha=(\alpha_1,\ldots, \alpha_d)$ and $\beta$ with $\alpha_i, \beta\in \NN_0$
 and $|\alpha|+2 \beta=m$ (where $|\alpha|=\sum_{i=1}^d\alpha_i$), 
\begin{align}\label{regu-est}
\norm{\nabla_x^{\alpha}\pd_t^{\beta}v}_{L^p(\R^d_+\times \R)}
\lesssim 
\norm{g}_{\dotc B^{m-\frac 1p, \frac m 2 - \frac 1{2p}}_{p,p} (\R^{d-1}\times \R)}.
\end{align}

\item For
$1\leq p\leq \infty $ and any $m\in \NN_0$,
\begin{align}\label{trace-est}
\norm{g}_{\dotc B^{m-\frac 1p, \frac m 2 - \frac 1{2p}}_{p,p} (\R^{d-1}\times \R)}\lesssim \norm{\pd_{x_d}^{m}v}_{L^p(\R^d_+\times \R)}.
\end{align}
\end{enumerate}
\end{thm}

Inequality \eqref{regu-est} is the regularity estimate of solution of heat equation in term of its boundary data, while inequality \eqref{trace-est} is the trace estimate which shows the optimality.
Note that we allow $p=1$ in both \eqref{regu-est} and \eqref{trace-est}, and we also allow $p=\infty$ in \eqref{trace-est}.
When $d=1$, $\dotc B^{s,s/2}_{p,p}(\R^{d-1}\times \R) = \dotc B^{s/2}_{p,p}(\R)$.

For corresponding estimates in the {time-independent} case, see Mironescu-Russ \cite[Theorem 1.9]{Mironescu-Russ-2015}, and Leoni \cite[Theorems 18.53, 18.57]{Leoni-2017}.
The related $W^{2,1}_{p,\rho}$ parabolic regularity estimate of $v$ in terms of its boundary data $g$ was first studied by
Weidemaier \cite{Weidemaier-02,Weidemaier-05} and Denk-Hieber-Pr\"uss \cite{DHP-07}, for $1<p,\rho<\infty$. For the borderline index $p,\rho=1$,  Ogawa-Shimizu \cite{Ogawa-Shimizu-2020-PJA}
and  \cite[Theorems 2.1 \& 2.2]{Ogawa-Shimizu-2022} established the following $L^1$-maximal regularity estimate for $1< p<\infty$, 
\begin{align}\label{shimizu-B}
	\norm{\pd_t v, \nabla^2 v}_{L^1(\R_+;\dotc B^0_{p,1}(\R^d_+))}\lesssim \norm{g}_{\dotc F^{1-\frac{1}{2p}}_{1,1}(\R_+;\dotc B^0_{p,1}(\R^{d-1}))}+\norm{g}_{L^1(\R_+;\dotc B^{2-\frac{1}{p}}_{p,1}(\R^{d-1}))},
\end{align}
and  for $1\le p<\infty$,
\begin{align}\label{shimizu1}
	\norm{\pd_t v, \nabla^2 v}_{L^1(\R_+;L^p(\R^d_+))}\lesssim \norm{g}_{\dotc F^{1-\frac{1}{2p}}_{1,1}(\R_+;L^p(\R^{d-1}))}+\norm{g}_{L^1(\R_+; B^{2-\frac{1}{p}}_{p,1}(\R^{d-1}))},
\end{align}
where $\dotc F^s_{p,q}$ is the Lizorkin-Triebel space.
Inequality \eqref{shimizu1} is similar to our regularity estimate \eqref{regu-est} when $m=2$ and $p=1$, and is different from \eqref{regu-est} in that \eqref{shimizu1} has a mixed exponent. 
Ogawa-Shimizu also established a trace estimate in \cite[Theorem 7.1]{Ogawa-Shimizu-2022}.
These results can be applied to prove the global well-posedness of free boundary problems for Navier-Stokes equations \cite{Ogawa-Shimizu-2024}.

With a different motivation, 
 a similar estimate of \eqref{shimizu1}  for the \emph{first order} derivative for $1\leq p<\infty$ was shown,\begin{equation}\label{clt2-1}
\norm{g}_{ \dotc{B}^{\frac 12-\frac 1{2p}}_{p,p}(\R)} \lesssim 
\norm{\partial_x v}_{L^p(\R_+ \times \R)} \lesssim \norm{g}_{ \dotc{B}^{\frac 12-\frac 1{2p}}_{p,p}(\R)}.
\end{equation}
The left inequality was first shown in Chang-Kang \cite[Eq.~(4.14)]{Chang-Kang-2023}. A different proof and the right inequality were given in Chen-Liang-Tsai \cite[Lemma 2.8]{Chen-Liang-Tsai-2025}. It is the key estimate to the construction of the derivative blow-up examples of Stokes system near the boundary.  
Inequalities \eqref{clt2-1} coincide with \eqref{regu-est} and \eqref{trace-est} when $m=1$ and $d=1$. Hence \cite[Lemma 2.8]{Chen-Liang-Tsai-2025} is a special case of Theorem \ref{thm_main}.

Unlike \cite{DHP-07,Ogawa-Shimizu-2020-PJA,Ogawa-Shimizu-2022,Weidemaier-02,Weidemaier-05}, our $g$ is defined for all $t\in\R$ and we have no initial data and no issue of compatibility at $t=0$.

A closer inspection of the proofs of \cite[Lemma 2.8]{Chen-Liang-Tsai-2025} and  \cite[Theorems 2.2 \& 7.1]{Ogawa-Shimizu-2022} reveals that they share the same spirit. We were thus motivated to find an estimate that encloses both \eqref{shimizu1} and \eqref{clt2-1}, but it has not been achieved yet. In this aspect, it is desired that
 $p=1$ is allowed in both \eqref{regu-est} and \eqref{trace-est}.

By applying the estimates in Theorem \ref{thm_main} and using interpolation, we get the following maximal regularity estimates of parabolic boundary value problems
in terms of homogeneous anisotropic Sobolev, Bessel potential and Besov norms (to be given in Section \ref{sec2}).

\begin{cor}\label{cor1}
Suppose $v$ satisfies \eqref{def-v}, which solves heat equation with Dirichlet data $g$.
\begin{enumerate}[(i)]
\item For $1\leq p< \infty$ and $m\in\NN_0$,
\begin{align}\label{est1}
	\norm{g}_{\dotc B^{2m-\frac 1p, m - \frac 1{2p}}_{p,p} (\R^{d-1}\times \R)}
 \lec 
\norm{v}_{\dotc W^{2m, m}_{p} (\R^d_+\times \R)}
 \lesssim 
 \norm{g}_{\dotc B^{2m-\frac 1p, m - \frac 1{2p}}_{p,p} (\R^{d-1}\times \R)},
\end{align}
where the left inequality also holds for $p=\infty$.

\item For $1< p< \infty$ and $s\geq 0$, 
\begin{align}\label{est2}
\norm{g}_{\dotc B^{s-\frac 1p, \frac s 2 - \frac 1{2p}}_{p,p} (\R^{d-1}\times \R)}
 \lec 
\norm{v}_{\dotc H^{s, \frac s 2}_{p} (\R^d_+\times \R)}
 \lesssim 
 \norm{g}_{\dotc B^{s-\frac 1p, \frac s 2 - \frac 1{2p}}_{p,p} (\R^{d-1}\times \R)}.
\end{align}

\item For $1\leq  p< \infty$ and $s> 0$, 
\begin{align}\label{est3}
\norm{g}_{\dotc B^{s-\frac 1p, \frac s 2 - \frac 1{2p}}_{p,p} (\R^{d-1}\times \R)}
 \lec 
\norm{v}_{\dotc B^{s, \frac s 2}_{p,p} (\R^d_+\times \R)}
 \lesssim 
 \norm{g}_{\dotc B^{s-\frac 1p, \frac s 2 - \frac 1{2p}}_{p,p} (\R^{d-1}\times \R)}.
\end{align}
\end{enumerate}
\end{cor}

The left inequalities in \eqref{est2} and \eqref{est3} are true for any $v\in \dotc H^{s, \frac s 2}_{p} (\R^d_+\times \R)$ and $v\in \dotc B^{s, \frac s 2}_{p,p} (\R^d_+\times \R)$ if $s>\frac{1}{p}$, $1<p<\infty$ and $g=v|_{x_d=0}$, without requiring $v$ to satisfy heat equation. See Amann \cite[Theorem 4.5.2]{Amann-2009}. We emphasize that \eqref{est1} and \eqref{est3} are valid for $p=1$ for a heat solution. 

Note that Bessel potential and Besov norms in some sense quantify the regularity the fractional derivatives, so from the heuristic of \cite{Chang-Kang-2023, Chen-Liang-Tsai-2025}, \eqref{est2} and \eqref{est3} might be useful for constructing  blow-up examples of Stokes system in terms of their fractional derivatives.

The preceding results can be easily extended to Neumann boundary condition $\pd_d v=g$ on $\pd \R^d_+\times \R$, in which case a solution is given by
\begin{align}\label{def-v2}
	v(x,t)=-2\int_\R\int_{\R^{d-1}}\Gamma(x'-y',x_d,t-s)g(y',s)\dy'\ds.
\end{align}
When $d=1$, the variable $y'$ is null and \eqref{def-v2} is reduced to $v(x,t)=-2\int_\R \Gamma(x,t-s)g(s)\ds$.

\begin{cor}\label{cor2}
Suppose $v$ satisfies \eqref{def-v2},  which solves heat equation with Neumann data $g$.
\begin{enumerate}[(i)]
\item For 
$1\leq p<\infty $, any $\alpha=(\alpha_1,\ldots, \alpha_d)$ and $\beta$ with $\alpha_i, \beta\in \NN_0$ and $|\alpha|+2 \beta=m$, 
\begin{align*}
\norm{\nabla_x^{\alpha}\pd_t^{\beta}v}_{L^p(\R^d_+\times \R)}
\lesssim 
\norm{g}_{\dotc B^{m-1-\frac 1p, \frac {m-1} 2 - \frac 1{2p}}_{p,p} (\R^{d-1}\times \R)}.
\end{align*}

\item For
$1\leq p\leq \infty $ and any $m\in \NN_0$
\begin{align*}
\norm{g}_{\dotc B^{m-1-\frac 1p, \frac {m-1} 2 - \frac 1{2p}}_{p,p} (\R^{d-1}\times \R)}\lesssim \norm{\pd_{x_d}^{m}v}_{L^p(\R^d_+\times \R)}.
\end{align*}

\item For $1\leq p< \infty$ and $m\in\NN_0$,
\begin{align*}
	\norm{g}_{\dotc B^{2m-1-\frac 1p, m-\frac {1} {2} - \frac 1{2p}}_{p,p} (\R^{d-1}\times \R)}
 \lec 
\norm{v}_{\dotc W^{2m, m}_{p} (\R^d_+\times \R)}
 \lesssim 
 \norm{g}_{\dotc B^{2m-1-\frac 1 p,m-\frac{1}{2} - \frac 1{2p}}_{p,p} (\R^{d-1}\times \R)},
\end{align*}
where the left inequality also holds for $p=\infty$.

\item For $1< p< \infty$ and $s\geq 0$, 
\begin{align*}
\norm{g}_{\dotc B^{s-1-\frac 1p, \frac {s-1} 2 - \frac 1{2p}}_{p,p} (\R^{d-1}\times \R)}
 \lec 
\norm{v}_{\dotc H^{s, \frac s 2}_{p} (\R^d_+\times \R)}
 \lesssim 
 \norm{g}_{\dotc B^{s-1-\frac 1p, \frac {s-1} 2 - \frac 1{2p}}_{p,p} (\R^{d-1}\times \R)}.
\end{align*}

\item For $1\leq  p< \infty$ and $s> 0$, 
\begin{align*}
\norm{g}_{\dotc B^{s-1-\frac 1p, \frac {s-1} 2 - \frac 1{2p}}_{p,p} (\R^{d-1}\times \R)}
 \lec 
\norm{v}_{\dotc B^{s, \frac s 2}_{p,p} (\R^d_+\times \R)}
 \lesssim 
 \norm{g}_{\dotc B^{s-1-\frac 1p, \frac {s-1} 2 - \frac 1{2p}}_{p,p} (\R^{d-1}\times \R)}.
\end{align*}
\end{enumerate}
\end{cor}

The rest of this note is structured as follows.
In Section \ref{sec2} we recall the definitions and properties (in particular their interpolation) of homogeneous
Sobolev, Bessel potential and Besov spaces, and their anisotropic versions.
In Section \ref{sec3} we prove the main theorem and the corollaries.

\section{Function spaces and interpolation}\label{sec2}
We first recall definitions and properties of classical homogeneous isotropic Sobolev, Bessel potential and Besov spaces in \S\ref{sec2.1}. 
We then recall their anisotropic versions  in \S\ref{sec2.2}. 

\subsection{Isotropic case}\label{sec2.1}
We first recall definitions and properties of the isotropic (without $t$ variable) function spaces. The standard reference is Triebel \cite{Triebel-1983}. 

In this subsection we denote by $\mathcal{F}$ the Fourier transform in $\R^d$, and denote by $\Delta_j$, $j \in \mathbb Z$, the usual Littlewood-Paley operators in $\R^d$. 
Denote by $\mathcal{S}'_h(\R^d)$ (\cite[Definition 1.26]{Bahouri-Chemin-Danchin-2011}) the subspace of tempered distributions $\mathcal{S}'(\R^d)$ consisting of those $u \in \mathcal{S}'(\R^d)$ s.t.
\[
	\lim_{\lambda \to \infty}\norm{\mathcal{F}^{-1}\Big(f(\lambda \xi) \cdot \cF u\Big)}_{L^\infty}=0,
\quad \text{for any } \ f\in C_c^\infty(\R^d).
\]
This space excludes some growth of the functions and makes the following homogeneous spaces (endowed with semi-norm) normed spaces.

\begin{defn}[Isotropic function spaces]\label{def-iso}\mbox{}
\begin{enumerate}[(i)]
\item For $1\leq p\leq \infty$ and $m\in \NN_0$, the \textit{homogeneous Sobolev space} $\dotc W^{m}_p(\R^d)$ consists of those distributions $u\in \mathcal{S}'_h(\R^d)$ s.t.
\begin{equation*}
	\norm{u}_{\dotc W^m_p(\R^d)}=\tsum_{|\alpha|=m}\norm{\nabla^{\alpha}u}_{L^p(\R^d)}<\infty.
\end{equation*}

\item For $1\leq p\leq  \infty$ and $s\in \R$, the \textit{homogeneous Bessel potential space} $\dotc H^s_p(\R^d)$ consists of those distributions $u\in \mathcal{S}'_h(\R^d)$ s.t.
\begin{align*}
	\norm{u}_{\dotc H^s_p(\R^d)}=\norm{\tsum_{j\in \mathbb{Z}}\mathcal{F}^{-1}\Big(|\xi|^s \cF (\Delta_j u)\Big)}_{L^p(\R^d)}<\infty.
\end{align*}

\item For $1\leq p,q \leq  \infty$ and $s\in \R$, the \textit{homogeneous Besov space} $\dotc B^s_{p,q}(\R^d)$ consists of those distributions $u\in \mathcal{S}'_h(\R^d)$ s.t.
\begin{align*}
	\norm{u}_{\dotc B^s_{p,q}(\R^d)}=\bigg(\tsum_{j\in \mathbb{Z}}2^{js q}\|\Delta_j u\|^q_{L^p(\R^d)}\bigg)^{1/q}<\infty,
\end{align*}
where for $q=\infty$ we take the norm to be $\sup_{j\in \mathbb{Z}} 2^{js}\|\Delta_j u\|_{L^p(\R^d)}$.


\item
For a smooth (bounded or unbounded) domain $\Omega\subset \R^d$ and 
\[X\in \{\dotc W^m_p\,,\dotc H^s_p\,, \dotc B^s_{p,q}\},\]
denote by $X(\Omega)$ the function space of those $u$ with finite norm,
\begin{align*}
	\norm{u}_{X(\Omega)}=\inf\{\norm{\tilde{u}}_{X(\R^d)}|\,\tilde{u}\in X(\R^d),\,\tilde{u}|_\Omega=u\}.
\end{align*}
\end{enumerate}
\end{defn}

Next, we recall some properties of these spaces. 

\begin{lem}[Basic properties]\label{fun-prop1}\mbox{}
\begin{enumerate}[(i)]
\item For $1<p<\infty$, $m\in \NN_0$ and $s\in \R$, 
\begin{align*}
\dotc W^m_p(\Omega)=\dotc H^m_p(\Omega),\quad \dotc H^s_2(\Omega)=\dotc B^s_{2,2}(\Omega).	
\end{align*} 

\item For $1\leq p\leq \infty$ and $m\in \NN_0$, $\dotc W^m_p(\R^d_+)$ has an equivalent norm
\begin{align*}
	\frac{1}{C} \norm{u}_{\dotc W^m_p(\R^d_+)}\leq\sum_{|\alpha|=m}\norm{\nabla^{\alpha}u}_{L^p(\R^d_+)}\leq C\norm{u}_{\dotc W^m_p(\R^d_+)}.
\end{align*}

\item For $1\leq p< \infty$ and $0<s\neq \hbox{integer}$, $\dotc B^s_{p,p}(\R^d_+)$ has an equivalent norm
\begin{align}\label{besov-explicit}
	\frac{1}{C} \norm{u}_{\dotc B^s_{p,p}(\R^d_+)}\leq\sum_{|\alpha|=\lfloor s\rfloor}\left(\int_{\R^d_+}\int_{\R^d_+}\frac{|\nabla^{\alpha}u(x)-\nabla^{\alpha}u(y)|^{p}}{|x-y|^{d+(s-\lfloor s\rfloor) p}}\dx\dy\right)^{1/p}\leq C\norm{u}_{\dotc B^s_{p,p}(\R^d_+)}.
\end{align}
\end{enumerate}
\end{lem}
\begin{proof}
See \cite[2.2.2, 2.3.1, 2.3.5, 3.4.2]{Triebel-1983} for related results concerning the inhomogeneous case. It's easy to see (i) is true for domain $\R^d$, and hence is true for $\Omega$. For (ii) and (iii), it's easy to see it's true for domain $\R^d$, and by high-order extension from $\R^d_+$ to $\R^d$ (which keeps the $m$-th derivative on the boundary $\pd \R^d_+$), it's also true in $\R^d_+$.
\end{proof}
\begin{rem}\label{rem2.3}
The domain $\R^d_+$ in Lemma \ref{fun-prop1} (ii) cannot be replaced by a bounded domain $\Omega$, since $u=\hbox{const}\neq 0$ has nontrivial $\dotc W^m_p(\Omega)$ norm. It's similar for Lemma \ref{fun-prop1} (iii).
See \cite[Proposition 3.4.2]{Triebel-1983} for the explicit nonhomogeneous norm $B^s_{p,p}(\Omega)$ in bounded domain. See \cite[Theorem 2.37]{Bahouri-Chemin-Danchin-2011} for the corresponding form of \eqref{besov-explicit} when $s$ is an integer.
\end{rem} 

Denote by $(\cdot,\cdot)_{\theta,q}$ the real interpolation and by $[\cdot,\cdot]_{\theta}$ the complex interpolation. In the following we assume $0<\theta<1$, and let
\begin{equation}\label{spq.def}
s=(1-\theta) s_0+\theta s_1,\quad \frac{1}{p}=\frac{1-\theta}{p_0}+\frac{\theta}{p_1},\quad \frac{1}{q}=\frac{1-\theta}{q_0}+\frac{\theta}{q_1}.
\end{equation}

\begin{lem}\label{fun-prop-interpo1} Let $\Om$ be a smooth \textup{(}bounded or unbounded\textup{)} domain in $\R^d$.
\begin{enumerate}[(i)]
\item For $1\leq p_0,p_1,q_0,q_1\leq \infty$ and $s_0,s_1\in \R$, 
\begin{align*}
[\dotc B^{s_0}_{p_0,q_0}(\Omega),\,\dotc B^{s_1}_{p_1,q_1}(\Omega)]_\theta=\dotc B^s_{p,q}(\Omega).
\end{align*}

\item For $1< p_0,p_1< \infty$ and $s_0,s_1\in \R$, 
\begin{align*}
[\dotc H^{s_0}_{p_0}(\Omega),\,\dotc H^{s_1}_{p_1}(\Omega)]_\theta=\dotc H^s_{p}(\Omega).
\end{align*}

\item For $1\leq p,q,q_0,q_1\leq \infty$ and $s_0\neq s_1\in \R$,
\begin{align*}
(\dotc B^{s_0}_{p,q_0}(\Omega),\,\dotc B^{s_1}_{p,q_1}(\Omega))_{\theta,q}=\dotc B^s_{p,q}(\Omega),
\end{align*}
\begin{align*}
(\dotc H^{s_0}_{p}(\Omega),\,\dotc H^{s_1}_{p}(\Omega))_{\theta,q}=\dotc B^s_{p,q}(\Omega).	
\end{align*} 

\item  For $1\leq p,q\leq \infty$ and $m=\NN_0\setminus \{0\}$,
\begin{align*}
(L^p(\Omega),\,\dotc W^m_{p}(\Omega))_{\theta,q}=\dotc B^{\theta m}_{p,q}(\Omega).
\end{align*}
\end{enumerate}
\end{lem}
\begin{proof}
(i),(ii) and (iii) follow directly from the inhomogeneous version in \cite[\S3.3.6]{Triebel-1983} (note $\dotc H^s_p=\dotc F^s_{p,2}$ in \cite[\S2.5.6]{Triebel-1983}). See \cite[\S 6.3, 6.4]{Bergh-Lofstrom-1976} for $p=1,\infty$ cases of $(iii)_2$. (iv) follows from \cite[Theorem 17.30]{Leoni-2017}.
\end{proof}

Bessel potential space $\dotc H^{s}_{p}$ and Besov space $\dotc B^{s}_{p,p}$ are two ways to quantify the properties of the fractional derivatives. Sometimes they are used to define the so-call \textit{fractional Sobolev space}. The advantage of Bessel potential space is that it matches with the usual Sobolev space when $s$ is an integer. The advantage of Besov space is that it has an explicit form \eqref{besov-explicit} of the norm (and \cite[Proposition 3.4.2]{Triebel-1983} for $B^{s}_{p,p}$) when $s$ is a non-integer real number.

\subsection{Anisotropic case}\label{sec2.2}

To deal with parabolic equations with time variable $t$, we introduce anisotropic function spaces parallel to those considered in the previous subsection. The standard reference is Amann \cite{Amann-2009}. See \cite[Section 2]{Chang-Kang-2018} for what are relevant to us.

Denote by $\cF_{x,t}$ the Fourier transform in $\R^d\times \R$. For complex number $a+ib$ and real number $s$, let $(a+ib)^{s}=e^{s\ln(a+ib)}$, where the branch cut of the logarithm is the negative real line. Define the anisotropic Littlewood-Paley operators in $\R^d\times \R$ as 
$\Delta_j u=\mathcal{F}_{x,t}^{-1}(\phi_j \cF_{x,t}(u))$, $j \in \mathbb Z$, 
with $\phi_j(\xi,\tau)= \phi(2^{-j}\xi,2^{-2j}\tau)$, and $\phi\in C_c^{\infty}(\R^{d}\times \R)$ satisfying
\begin{align*}
    \left\{
\begin{aligned}
\phi(\xi, \tau)>0, \quad &\hbox{on }2^{-1}<|\xi|+|\tau|^{1/2}<2, \\
\phi(\xi, \tau)=0, \quad &\hbox{elsewhere,}\\
\sum_{j\in \mathbb{Z}} \phi_j(\xi,\tau)=1,\quad &\hbox{on }(\xi,\tau)\neq (0,0).
\end{aligned} \right. 
\end{align*}

\begin{defn}[Anisotropic function spaces]\label{def-aniso}\mbox{}
\begin{enumerate}[(i)]
\item For $1\leq p\leq \infty$ and $m\in \NN_0$, the \textit{anisotropic homogeneous Sobolev space} $\dotc W^{2m,m}_p(\R^d\times \R)$ consists of those distributions $u\in \mathcal{S}'_h(\R^d\times \R)$ s.t.
\begin{align*}
	\norm{u}_{\dotc W^{2m,m}_p(\R^d\times \R)}=\sum_{|\alpha|+2\beta=2m}\norm{\nabla_x^{\alpha}\pd_t^{\beta}u}_{L^p(\R^d\times \R)}<\infty.
\end{align*}

\item For $1\leq p\leq  \infty$ and $s\in \R$, the \textit{anisotropic homogeneous Bessel potential space} $\dotc H^{s,\frac{s}{2}}_p(\R^d\times \R)$ consists of those distributions $u\in \mathcal{S}'_h(\R^d\times \R)$ s.t.
\begin{align*}
	\norm{u}_{\dotc H^{s,\frac{s}{2}}_p(\R^d\times \R)}=\norm{\tsum_{j\in\mathbb{Z}}\mathcal{F}_{x,t}^{-1}\Big((|\xi|^2+i\tau)^{\frac{s}{2}}\,\cF_{x,t} (\Delta_ju)\Big)}_{L^p(\R^d\times \R)}<\infty.
\end{align*}

\item For $1\leq p,q \leq  \infty$ and $s\in \R$, the \textit{anisotropic homogeneous Besov space} $\dotc B^{s,\frac{s}{2}}_{p,q}(\R^d\times \R)$ consists of those distributions $u\in \mathcal{S}'_h(\R^d\times \R)$ s.t.
\begin{align*}
	\norm{u}_{\dotc B^{s,\frac{s}{2}}_{p,q}(\R^d\times\R)}=\left(\tsum_{j\in \mathbb{Z}}2^{js q}\|\Delta_j u\|^q_{L^p(\R^d\times \R)}\right)^{1/q}<\infty,
\end{align*}
where for $q=\infty$ we take the norm to be $\sup_{j\in \mathbb{Z}} 2^{js}\|\Delta_j u\|_{L^p(\R^d\times\R)}$.

\item
For $X\in \{\dotc W^{2m,m}_p\,,\dotc H^{s,\frac{s}{2}}_p\,, \dotc B^{s,\frac{s}{2}}_{p,q}\}$,
denote by $X(\R^d_+\times \R)$ the function space of those $u$ s.t. 
\begin{align*}
	\norm{u}_{X(\R^d_+\times \R)}=\inf\{\norm{\tilde{u}}_{X(\R^d\times \R)}|\,\tilde{u}\in X(\R^d\times \R),\,\tilde{u}|_{\R^d_+\times \R}=u\}<\infty.
\end{align*}
\end{enumerate}
\end{defn}

Next we list some properties of the above spaces.

\begin{lem}\label{fun-prop2}
Let $\mathbb Y=\R^{d}$ or $\mathbb Y=\R^d_+$, $d \ge1$.
Assume $0<\theta<1$ and $s,p,q$ satisfy \eqref{spq.def}.

\begin{enumerate}[(i)]
\item For $1<p<\infty$ and $m\in \NN_0$, 
\begin{align*}
\dotc W^{2m,m}_p(\mathbb Y\times \R)=\dotc H^{2m,m}_p(\mathbb Y\times \R).
\end{align*}

\item For $1< p,q< \infty$ and $s_0\neq s_1\in \R$, 
\begin{align*}
\Big[\dotc B^{s_0,\frac{s_0}{2}}_{p,q}(\mathbb Y\times \R),\,\dotc B^{s_1,\frac{s_1}{2}}_{p,q}(\mathbb Y\times \R)\Big]_\theta=\dotc B^{s,\frac{s}{2}}_{p,q}(\mathbb Y\times \R).	
\end{align*}

\item For $1< p_0,p_1< \infty$ and $s_0\neq s_1\in \R$,
\begin{align*}
\Big[\dotc H^{s_0,\frac{s_0}{2}}_{p_0}(\mathbb Y\times \R),\,\dotc H^{s_1,\frac{s_1}{2}}_{p_1}(\mathbb Y\times \R)\Big]_\theta=\dotc H^{s,\frac{s}{2}}_{p}(\mathbb Y\times \R).	
\end{align*}

\item For $1< p,q,q_0,q_1< \infty$ and $s_0\neq s_1\in \R$,
\begin{align*}
\Big(\dotc B^{s_0,\frac{s_0}{2}}_{p,q_0}(\mathbb Y\times \R),\,\dotc B^{s_1,\frac{s_1}{2}}_{p,q_1}(\mathbb Y\times \R)\Big)_{\theta,q}&=\dotc B^{s,\frac{s}{2}}_{p,q}(\mathbb Y\times \R),
\\[4pt]
\Big(\dotc H^{s_0,\frac{s_0}{2}}_{p}(\mathbb Y\times \R),\,\dotc H^{s_1,\frac{s_1}{2}}_{p}(\mathbb Y\times \R)\Big)_{\theta,q}&=\dotc B^{s,\frac{s}{2}}_{p,q}(\mathbb Y\times \R).	
\end{align*} 

\item  Let
\begin{align*}
\dotc W^{2m,m}_{p}(\mathbb Y\times \R)\subset \dotc W^{2m,m,*}_{p}(\mathbb Y\times \R)\equiv L^p(\R, \dotc W^{2m}_{p}(\mathbb Y))\cap L^p(\mathbb Y, \dotc W^{m}_{p}(\R)).
\end{align*}
Then we have, for $1\leq p<\infty$ and $m\in \NN_0\setminus \{0\}$,
\begin{align*}
\Big(L^p(\mathbb Y\times \R),\,\dotc W^{2m,m,*}_{p}(\mathbb Y\times \R)\Big)_{\theta,p}=\dotc B^{2 \theta m,\theta m}_{p,p}(\mathbb Y\times \R).
\end{align*}
\end{enumerate}
\end{lem}
\begin{proof}
(i) is well-known. For (ii), (iii) and (iv), see
\cite[(3.3.12), (3.4.1), Theorem 3.7.1]{Amann-2009} for the case of the whole space,
and \cite[Theorem 4.4.1]{Amann-2009} for the half-space case. 

For (v), by Lemma \ref{fun-prop-interpo1} (iv) and \cite[\S5.8.6]{Bergh-Lofstrom-1976} we have for $m\in \NN_0\setminus \{0\}$,
\begin{align*}
\Big(L^p(\R, L^p(\mathbb Y)),\,L^p(\R, \dotc W^{2m}_{p}(\mathbb Y))\Big)_{\theta,p}&=L^p(\R, \dotc B^{2 \theta m}_{p,p}(\mathbb Y)),\\[4pt]
\Big(L^p(\mathbb Y, L^p(\R)),\,L^p(\mathbb Y, \dotc W^{m}_{p}(\R))\Big)_{\theta,p}&=L^p(\mathbb Y, \dotc B^{\theta m}_{p,p}(\R)).
\end{align*}
By \cite[Theorem 3.6.3]{Amann-2009} we have for $1\leq p<\infty$,
\begin{align*}
	\dotc B^{2 \theta m,\theta m}_{p,p}(\mathbb Y\times \R)=L^p(\R, \dotc B^{2 \theta m}_{p,p}(\mathbb Y))\cap L^p(\mathbb Y, \dotc B^{\theta m}_{p,p}(\R)).
\end{align*}
Hence, we reach (v).
\end{proof}

\section{Maximal regularity estimate}\label{sec3}

In this section we prove Theorem \ref{thm_main}, Corollary \ref{cor1} and Corollary \ref{cor2}. Throughout this section we denote by $\hat{u}=\cF_{x',t}(u)$ the Fourier transform of $u$ in $(x',t)$ variables ($\hat{u}=\cF_{t}(u)$ when $d=1$), and $ (\hat u)^\vee=\cF_{x',t}^{-1}\, (\hat u)$. 
We also denote by $\Delta_j$ the anisotropic Littlewood-Paley operator in $(x',t)\in \R^{d-1}\times \R$, given by
$\Delta_j f  = ( \phi_j \hat{f})^\vee$ for $ \phi(\xi',\tau)$ satisfying  properties slightly different to those in Definition \ref{def-aniso} (iii), which is defined in $(x,t)\in \R^{d}\times \R$.

\begin{proof}[Proof of Theorem \ref{thm_main}]

The proof is similar to Chen-Liang-Tsai \cite[Lemma 2.8]{Chen-Liang-Tsai-2025}.  Recall that the Fourier transform of heat kernel (\cite[Corollary 2.5]{Chen-Liang-Tsai-2025}) for $x_d>0$ is
\[
\mathcal{F}_{x',t}(\pd_{x_d}\Gamma)=-\frac{1}{2(2\pi)^{d/2}}e^{-x_d\sqrt{|\xi'|^2+i\tau}}.
\]
When $d=1$, the variable $\xi'$ is null and $\mathcal{F}_{t}(\pd_{x_d}\Gamma)=\frac{-1}{2(2\pi)^{1/2}}e^{-x_d\sqrt{i\tau}}$.
Hence for $v$ in \eqref{def-v} we have that (note $\widehat{f*g}=(2\pi)^{d/2}\hat{f}\hat{g}$ by the definition of Fourier transform in \cite[\S 2]{Chen-Liang-Tsai-2025})
\begin{align}\label{v-g-hat}
    \hat{v}=  e^{-x_d\sqrt{|\xi'|^2+i\tau}} \cdot \hat{g}(\xi',\tau).
\end{align}

\noindent \textbf{\textit{Step 1}}.\quad Proof of \eqref{regu-est}.

\smallskip

For fixed $x_d>0$, we first claim that
\begin{align}\label{est-fix-xn}
\norm{\nabla_x^{\alpha}\pd_t^{\beta}v}_{L^p(\R^{d-1}\times \R)} \lesssim \tsum_{j\in \mathbb{Z}}\ 2^{m j}e^{- c  x_d2^{j}} \|\Delta_j  g\|_{L^{p}(\R^{d-1}\times \R)}.
\end{align}
By \eqref{v-g-hat} we have
\begin{equation*}
    \cF_{x',t}\Big(\nabla_x^{\alpha}\pd_t^{\beta}v\Big)=P(\xi',\tau)\cdot e^{-x_d\sqrt{|\xi'|^2+i\tau}} \cdot \hat{g},
\end{equation*}
with
\begin{equation*}
    P(\xi',\tau)\equiv
   \left( {\textstyle \prod}_{k=1}^{d-1}(i\xi_k)^{\alpha_k}\right)\cdot (-\sqrt{|\xi'|^2+i\tau})^{\alpha_d}\cdot (i\tau)^{\beta},  
\end{equation*}
which implies that
\begin{align*}
\|\Delta_j \nabla_x^{\alpha}\pd_t^{\beta}v\|_{L^{p}(\R^{d-1}\times \R)}=\|\Big({\phi}_j\cdot P(\xi',\tau)\cdot e^{-x_d\sqrt{|\xi'|^2+i\tau}}\cdot \hat{g}\Big)^{\vee}\|_{L^{p}(\R^{d-1}\times \R)}. 
\end{align*}
Let ${\psi}(\xi',\tau)$ be smooth and supported in some annulus with ${\psi}=1$ on the support of ${\phi}$.
 By Young's inequality and scaling technique, we can write the above equation as
\begin{align*}
&\,\quad \|\Delta_j \nabla_x^{\alpha}\pd_t^{\beta} v\|_{L^{p}(\R^{d-1}\times \R)}\\[4pt]
&=\|\Big({\psi}(2^{-j}\xi',2^{-2j}\tau)\cdot{\phi}(2^{-j}\xi',2^{-2j}\tau)\cdot P(\xi',\tau)\cdot e^{-x_d\sqrt{|\xi'|^2+i\tau}}\cdot \hat{g}\Big)^{\vee}\|_{L^{p}(\R^{d-1}\times \R)}\\[4pt]
&\leq \|\Big({\psi}(2^{-j}\xi',2^{-2j}\tau)\cdot P(\xi',\tau)\cdot e^{-x_d\sqrt{|\xi'|^2+i\tau}}\Big)^{\vee}\|_{L^{1}(\R^{d-1}\times \R)}\cdot \|\Delta_j  g\|_{L^p(\R^{d-1}\times \R)}\\[4pt]
&= 2^{m j}\, \|\Big({\psi}(\xi',\tau)\cdot P(\xi',\tau)\cdot e^{-x_d2^{j}\sqrt{|\xi'|^2+i\tau}}\Big)^{\vee}\|_{L^{1}(\R^{d-1}\times \R)}\cdot \|\Delta_j  g\|_{L^{p}(\R^{d-1}\times \R)}. 
\end{align*}
Therefore, what remains is to show that for any $x_d>0$,
\begin{align}\label{e-decay}
\|\Big({\psi}(\xi',\tau)\cdot P(\xi',\tau)\cdot e^{-x_d\sqrt{|\xi'|^2+i\tau}}\Big)^{\vee}\|_{L^{1}(\R^{d-1}\times \R)}\lesssim e^{- c  x_d}. 
\end{align}
Using integration by parts, we get
\begin{align}\label{a1}
&\,\quad \Big({\psi}(\xi',\tau)\cdot P(\xi',\tau)\cdot e^{-x_d\sqrt{|\xi'|^2+i\tau}}\Big)^{\vee}\notag\\[4pt]
&=\frac{1}{(2\pi)^{d/2}}\int_{\R^{d-1}\times \R}{\psi}\cdot P(\xi',\tau)\cdot e^{-x_d\sqrt{|\xi'|^2+i\tau}}\cdot e^{i(\xi'\cdot x'+\tau t )}\dxi'\dtau\notag\\[4pt]
&=\frac{1}{(2\pi)^{d/2}(1+\sum_{k=1}^{d-1}|x_k|^{2d}+t^{2d})}\int_{\hbox{supp}({\psi})}{\psi}\cdot P(\xi',\tau)\cdot e^{-x_d\sqrt{|\xi'|^2+i\tau}}\notag\\[4pt]
&\,\quad \times\Big(1+\sum_{k=1}^{d-1}(-\pd_{\xi_k}^2)^d+(-\pd_{\tau}^2)^d\Big)e^{i(\xi'\cdot x'+\tau t )}\dxi'\dtau\notag\\[4pt]
&=\frac{1}{(2\pi)^{d/2}(1+\sum_{k=1}^{d-1}|x_k|^{2d}+t^{2d})}\int_{\hbox{supp}({\psi})}e^{i(\xi'\cdot x'+\tau t )}\notag\\[4pt]
&\,\quad \times\Big(1+\sum_{k=1}^{d-1}(-\pd_{\xi_k}^2)^d+(-\pd_{\tau}^2)^d\Big)\Big({\psi}\cdot P(\xi',\tau)\cdot e^{-x_d\sqrt{|\xi'|^2+i\tau}}\Big)\dxi'\dtau. 
\end{align}
Note that $\hbox{Re}(\sqrt{|\xi'|^2+i\tau})\geq \frac{1}{\sqrt 2}(|\xi'|^4+|\tau|^2)^{1/4}$ no matter $\tau>0$ or $\tau<0$,  so we have
\begin{align*}
\Big|\Big(1+\sum_{k=1}^{d-1}(-\pd_{\xi_k}^2)^d+(-\pd_{\tau}^2)^d\Big)\Big({\psi}\cdot P(\xi',\tau)\cdot e^{-x_d\sqrt{|\xi'|^2+i\tau}}\Big)\Big|\lesssim e^{-cx_d}, 
\end{align*}
which together with \eqref{a1} leads to \eqref{e-decay}, and hence \eqref{est-fix-xn}. 

Now we are in position to show  \eqref{regu-est}. By virtue of \eqref{est-fix-xn}, we have
\begin{align*}
&\,\quad \int_{0}^{+\infty}\norm{\nabla_x^{\alpha}\pd_t^{\beta}v}_{L^p(\R^{d-1}\times \R)}^p \dx_d\\[4pt]
&\lesssim \int_{0}^{+\infty}\left(\tsum_{j\in \mathbb{Z}}\ 2^{m j}e^{- c  x_d2^{j}} \|\Delta_j  g\|_{L^{p}(\R^{d-1}\times \R)}\right)^p \dx_d  \\[4pt]
&\lesssim \tsum_{j\in \mathbb{Z}}\int_{0}^{+\infty} 2^{m jp -jp}\cdot x_d^{-1+1/p} \left(2^{ j}e^{-cx_d 2^{ j}} \right)^{p-1+1/p} \dx_d\cdot\|\Delta_j  g\|_{L^{p}(\R^{d-1}\times \R)}^p  \\[4pt]
&\lesssim \tsum_{j\in \mathbb{Z}}\ 2^{m jp-j}\ \|\Delta_j  g\|_{L^{p}(\R^{d-1}\times \R)}^p
= \|g\|_{\dotc{B}^{m-\frac 1p, \frac m 2 - \frac 1{2p}}_{p,p}}^p,
\end{align*}
which gives rise to \eqref{regu-est}. Note that the second inequality above is by H\"{o}lder's inequality when $p>1$,
\begin{align*}
	(\tsum a_jb_j)^p\leq (\tsum a_j^p)(\tsum b_j^{p/(p-1)})^{p-1}, 
\end{align*}
where 
\[a_j=2^{m j -j}\cdot \left(2^{j}e^{-cx_d 2^{j}} \right)^{(p-1+1/p)/p}\|\Delta_j  g\|_{L^{p}(\R^{d-1}\times \R)},
\quad b_j=\left(2^{j}e^{-cx_d 2^{j}} \right)^{(p-1)/p^2},
\]
and $\sum b_j^{p/(p-1)} \lesssim x_d^{-1/p}$ since $\sup_{x_d>0} \sum_{j\in \mathbb{Z}}(x_d2^{j}e^{-cx_d 2^{j}})^{1/p}\lesssim 1$  (see \cite[Lemma 2.35]{Bahouri-Chemin-Danchin-2011}). When $p=1$,  the second inequality above is trivial and can be skipped. 

\bigskip

\noindent \textbf{\textit{Step 2}}.\quad Proof of \eqref{trace-est}.

\smallskip

By \eqref{v-g-hat} we have
\[
\hat{g}(\xi',\tau)=4\int_{0}^{+\infty} x_d (|\xi'|^2+i\tau)\cdot e^{-x_d\sqrt{|\xi'|^2+i\tau}}\cdot \hat{v}(\xi',x_d,\tau) \dx_d. 
\]
Therefore, by Minkowski and Young's inequalities, \eqref{e-decay} and  the scaling argument above it, we have
\begin{align*}
&\,\quad \|\Delta_j g\|_{L^{p}(\R^{d-1}\times \R)}=\|\Big({\phi}_j\hat{g}\Big)^{\vee}\|_{L^{p}(\R^{d-1}\times \R)}\\[4pt]
&=\|4\int_0^{\infty}\Big({\phi}_j\cdot x_d (-\sqrt{|\xi'|^2+i\tau})^{2-m}\cdot e^{-x_d\sqrt{|\xi'|^2+i\tau}}\cdot(-\sqrt{|\xi'|^2+i\tau})^{m}\cdot \hat{v}\Big)^{\vee}\dx_d\|_{L^{p}(\R^{d-1}\times\R)}\\[4pt]
&\lesssim\int_0^{\infty}\norm{\big({\phi}_j\cdot x_d (-\sqrt{|\xi'|^2+i\tau})^{2-m}\cdot e^{-x_d\sqrt{|\xi'|^2+i\tau}}\big)^\vee }_{L^{1}(\R^{d-1}\times\R)}\cdot\|\pd_{x_d}^mv\|_{L^{p}(\R^{d-1}\times\R)}\dx_d   \\[4pt]
&\lesssim \int_{0}^{+\infty} x_d2^{(2-m)j} e^{-cx_d2^{j}}\|\pd_{x_d}^mv\|_{L^{p}(\R^{d-1}\times\R)}\dx_d.
\end{align*}
Note that in deriving \eqref{e-decay} we only used the scaling property of derivatives of $P(\xi',\tau)$, not its explicit form. In particular, $2-m$ is allowed to be negative.

Thus for $p=\infty$, we directly get 
\begin{align*}
\sup_{j\in\mathbb{Z}}2^{mj}\|\Delta_j g\|_{L^{\infty}(\R^{d-1}\times \R)}\lesssim \|\pd_{x_d}^mv\|_{L^{\infty}(\R^{d}_+\times\R)},
\end{align*}
and for $1\leq p<\infty$,
\begin{align*}
&\,\quad \tsum_{j\in \mathbb{Z}}\ 2^{(m p -1)j}\ \|\Delta_j  g\|_{L^{p}(\R^{d-1}\times \R)}^p \\[4pt]
&\lesssim \tsum_{j\in \mathbb{Z}}\ 2^{(m p -1)j}\ \left(\int_{0}^{+\infty} x_d2^{(2-m)j} e^{-cx_d2^{j}}\|\pd_{x_d}^mv\|_{L^{p}(\R^{d-1}\times\R)}\dx_d\right)^p  \\[4pt]
&\lesssim \int_{0}^{+\infty}\left(\tsum_{j\in \mathbb{Z}}\ x_d2^{j}\,e^{-cx_d2^{j}}\right)\cdot   \|\pd_{x_d}^mv\|_{L^{p}(\R^{d-1}\times\R)}^p\dx_d\\[4pt]
&\lesssim \int_{0}^{+\infty}   \|\pd_{x_d}^mv\|_{L^{p}(\R^{d-1}\times\R)}^p\dx_d,
\end{align*}
which gives rise to \eqref{trace-est}. For the second inequality above we have used H\"{o}lder's inequality 
\[
(\int a(x)b(x)\dx)^p\leq (\int a(x)^p\dx)(\int b(x)^{p/(p-1)}\dx)^{p-1},
\] 
where 
\[a=\big(x_d2^{j}\,e^{-cx_d2^{j}}\big)^{1/p}\|\pd_{x_d}^mv\|_{L^{p}(\R^{d-1}\times\R)},\quad
b=2^{(1-m)j}\big(x_d2^{j}\,e^{-cx_d2^{j}}\big)^{1-1/p}.
\]
For the third inequality we used $\sup_{x>0} \sum_{j\in \mathbb{Z}}x2^{j}e^{-cx 2^{j}}\lesssim 1$ again.
\end{proof}

\begin{proof}[Proof of Corollary \ref{cor1}]
It's easy to see \eqref{est1} follows from Theorem \ref{thm_main}. Estimates \eqref{est2} follows from \eqref{est1} and Lemma \ref{fun-prop2} (i)--(iii). To prove estimates \eqref{est3}, firstly it's easy to see \eqref{est1} is also true if we replace $\dotc W^{2m,m}_{p}(\R^d_+\times \R)$ with $\dotc W^{2m,m,*}_{p}(\R^d_+\times \R)$, defined in Lemma \ref{fun-prop2} (v). Next we apply Lemma \ref{fun-prop2} (v). When $d=1$, $\dotc B^{s,s/2}_{p,p}(\R^{d-1}\times \R) = \dotc B^{s/2}_{p,p}(\R)$ and we also use Lemma \ref{fun-prop-interpo1}.
\end{proof}

\begin{proof}[Proof of Corollary \ref{cor2}]
For $v$ satisfying \eqref{def-v2}, we have (see \cite[Lemma 2.4]{Chen-Liang-Tsai-2025}) 
\begin{align*}
    \hat{v}=  -\frac{1}{\sqrt{|\xi'|^2+i\tau}}e^{-x_d\sqrt{|\xi'|^2+i\tau}} \cdot \hat{g}(\xi',\tau). 
\end{align*}
The rest of the proof is the same as Theorem \ref{thm_main} and Corollary \ref{cor1}.
\end{proof}

\section*{Acknowledgments}
We warmly thank Professor Senjo Shimizu for fruitful discussions. TT also thanks colleagues in Kyoto University for their kind hospitality during his visit in June 2025.

Chen was supported in part by National Natural Science Foundation of China under grant [12101556] and Zhejiang Provincial Natural Science Foundation of China under  grant [LY24A010015] and Fundamental Research Funds of ZUST [2025QN062]. 
The research of both Liang and Tsai was partially supported by Natural Sciences and Engineering Research Council of Canada (NSERC) under grant RGPIN-2023-04534.
\bibliography{MaxReg}
\bibliographystyle{abbrv}

\end{document}